\documentclass{mcom-l}

\newtheorem{theorem}{Theorem}[section]
\newtheorem{lemma}[theorem]{Lemma}

\theoremstyle{definition}

\theoremstyle{remark}

\numberwithin{equation}{section}

\newcommand{\nvec}{\ensuremath{\mathbf{n}}}
\newcommand{\uvec}{\ensuremath{\mathbf{u}}}
\newcommand{\vvec}{\ensuremath{\mathbf{v}}}
\newcommand{\wvec}{\ensuremath{\mathbf{w}}}
\newcommand{\xvec}{\ensuremath{\mathbf{x}}}

\DeclareMathOperator{\Div}{div}

\newcommand{\Bh}{\ensuremath{\mathcal{B}^3(\mathcal{T}_h)}}
\newcommand{\Eh}{\ensuremath{\mathcal{E}_h}}
\newcommand{\Nh}{\ensuremath{\mathcal{N}_h}}
\newcommand{\Sh}{\ensuremath{\mathcal{S}_h}}
\newcommand{\Th}{\ensuremath{\mathcal{T}_h}}

\newcommand{\R}{\ensuremath{\mathbb{R}}}

\begin{document}

\title[Mixed finite element discretization for nonlinear problems]
      {Efficient Realization of the Mixed Finite Element Discretization 
       for nonlinear Problems}

\author{Peter Knabner}
\address{Institute for Applied Mathematics, University Erlangen-N\"urnberg,
Martensstr.~3, D-91058 Erlangen, Germany}
\email{knabner@am.uni-erlangen.de}

\author{Gerhard Summ}
\address{Institute for Applied Mathematics, University Erlangen-N\"urnberg,
Martensstr.~3, D-91058 Erlangen, Germany}
\email{summ@am.uni-erlangen.de}

%    General info

\date{\today}

%\dedicatory{This paper is dedicated to our authors.}


\begin{abstract}
We consider implementational aspects of the mixed finite element method
for a special class of nonlinear problems.
We establish the equivalence of the hybridized formulation of the mixed
finite element method to a nonconforming finite element method with
augmented Crouzeix-Raviart ansatz space.
We discuss the reduction of unknowns by static condensation
and propose Newton's method for the solution of local and global systems.
Finally, we show, how such a nonlinear problem arises 
from the mixed formulation of Darcy--Forchheimer flow in porous media.
\end{abstract}

\maketitle

\section{Introduction}
The mixed finite element method is widely used
for the discretization of flow problems,
since it provides a direct and accurate approximation of the flux.
However, the mixed finite element discretization gives rise to large
systems of algebraic equations, which are difficult to solve,
because they are of saddle point type.
For linear problems the introduction of inter-element multipliers,
commonly called hybridization, is frequently used to transform this saddle
point problem into a problem, whose matrix is symmetric and positive definite
(cf.\ \cite[\S~V.1]{Brezzi/Fortin}).
Furthermore static condensation can be used to reduce the size of the system,
which finally has to be solved.
Exploiting the equivalence of the hybridized mixed finite element method
to certain nonconforming finite element methods,
several authors develop multigrid methods to solve the resulting system 
of linear equations (see e.g.\ \cite{Arbogast/Chen:95}, 
\cite{Brenner:92}, \cite{Chen:96}).
As far as we know, for nonlinear problems results of this kind 
are still missing.
In this article, we describe,
how the above mentioned techniques can be generalized to a certain class
of doubly nonlinear equations of monotone type.
We consider a bounded domain $\Omega$ with polygonal boundary $\partial \Omega$.
As usual, the Lebesgue spaces are denoted by $L^r(\Omega)$ 
for $1 \le r \le \infty$ and equpped with the norm $\| \cdot \|_{0,r,\Omega}$.
$\left( L^s(\Omega) \right)^2$ is the space of vector functions $\vvec$,
whose components belong to $L^s(\Omega)$.

We consider the variational problem:\\
Find $(\uvec,p) \in V \times Q$ such that
\begin{subequations}		\label{gen.var.pr}
\begin{align}			\label{gen.var.pr:a}
a(\uvec,\vvec) - b(\vvec,p) & = g(\vvec) \quad \text{for all } \vvec \in V \: ,\\
c(p,q) + b(\uvec,q) & = f(q) \quad \: \text{for all } q \in Q \: . \label{gen.var.pr:b}
\end{align}
\end{subequations}
Here $Q := L^r(\Omega)$ and $V$ is defined by
\begin{equation}			\label{def:V}
V := \left\{ \vvec \in \left( L^s(\Omega) \right)^n \bigm|
	     \Div \vvec \in L^{r'}(\Omega) \right\} \: ,
\end{equation}
where $r,s \in (1,\infty)$ and $1/r+1/r'=1$.
The bilinear form $b$ on $V \times Q$ is 
defined by $b(\vvec,q) := \int_\Omega \Div \vvec q \, d\xvec$,
$g$ and $f$ are linear forms on $V$ and $Q$, resp., and
$a$ and $c$ are continuous (nonlinear) forms on $V \times V$ and $Q \times Q$,
resp., defined by
$$ a(\uvec,\vvec) := \int_\Omega G(\uvec) \cdot \vvec \, d\xvec \quad \text{and}
   \quad c(p,q) := \int_\Omega R(p) \, q \, d\xvec \: , $$
where $R: L^r(\Omega) \to L^{r'}(\Omega)$ is a continuous mapping and
$G : \left( L^s(\Omega) \right)^n \to \big( L^{s'}(\Omega) \big)^n$
is a strictly monotone, coercive and continuous mapping.
A variational problem like (\ref{gen.var.pr}) occurs as mixed formulation
\begin{equation}			\label{mixed.form}
 \left. \begin{array}{r@{\:=\:}l}
    G(\uvec) + \nabla p & 0 \\
    R(p) + \Div \uvec & f
  \end{array} \right\} \quad \text{in } \Omega
\end{equation}
of the doubly nonlinear elliptic equation
\begin{equation}			\label{ell.eq} 
R(p) - \Div G^{-1}(\nabla p) = f \quad \text{in } \Omega
\end{equation}
with Dirichlet boundary conditions $p=g$ on $\partial\Omega$.
Then the linear forms $g$ and $f$ are defined by
$$ g(\vvec) := \int_{\partial\Omega} g \, (\vvec \cdot \nvec) \, d\sigma \quad
   \text{and} \quad f(q) := \int_\Omega f q \, d\xvec \: . $$

%%%%%%%%%%%%%%%%%%%%%%%%%%%%%%%%%%%%%%%%%%%%%%%%%%%%%%%%%%%%%%%%%%%%%%
%
\section{Discretization with the mixed finite element method}
%
%%%%%%%%%%%%%%%%%%%%%%%%%%%%%%%%%%%%%%%%%%%%%%%%%%%%%%%%%%%%%%%%%%%%%%
Replacing the infinite-dimensional spaces $V$ and $Q$ in (\ref{gen.var.pr})
by finite-dimensional subspaces $V_h \subset V$ and $Q_h \subset Q$ we
arrive at the mixed finite element method.

%%%%%%%%%%%%%%%%%%%%%%%%%%%%%%%%%%%%%%%%%%%%%%%%%%%%%%%%%%%%%%%%%%%%%%
%
\subsection{The finite-dimensional spaces $V_h$ and $Q_h$}
The solution of nonlinear problems like (\ref{ell.eq}) is, 
in general, not very regular.
In addition, we have to deal with nonsmooth parameter functions in our 
applications. 
Thus an approximation by higher-order finite elements makes no sense.
Therefore we choose the Raviart--Thomas space of lowest order, i.e.
\begin{subequations}				\label{disc_spaces}
\begin{align}						\label{V_h}
V_h & := RT_0(\Omega,\Th) =
	   \left\{ \vvec_h \in V \bigm|
	           \vvec_h |_K \in RT_0(K) 
	           \text{ for all } K \in \Th \right\} \, , \\
Q_h & := \left\{ q_h \in Q \bigm| q_h |_K \in P_0(K)	
	           \text{ for all } K \in \Th \right\} \, , 	\label{Q_h}
\end{align}
\end{subequations}
where $RT_0(K)$ is defined by
$$ RT_0(K) := \big( P_0(K) \big)^2 + 
   \Big( \! \begin{array}{c} x \\[-0.8ex] y \end{array} \! \Big) P_0(K) \: . $$
Here $P_k(K)$ denotes the space of polynomials of degree $\le k$ on $K$;
hence $P_0(K)$ is the space of constants on $K$.
The definition of $V_h$ and $Q_h$ is based on a triangulation $\Th$ 
of $\Omega$ into triangular elements $K$.
We denote by $\Eh$ the set of edges of $\Th$, which can be subdivided into
the set of interior edges
$ \Eh^i := \{ e \in \Eh \, | \, e \not \subset \partial \Omega \} $
and the set of boundary edges
$ \Eh^b := \{ e \in \Eh \, | \, e \subset \partial \Omega \}\,$.
In addition, we define a unit normal vector $\nvec_e$ for every $e \in \Eh$.
The orientation of $\nvec_e$ is arbitrary.

The requirement $\vvec_h \in V$ implies that the normal components
of the fluxes across interior edges are continuous:

\begin{lemma}				\label{lemma:W^3(div):char}
The following definition for $V$ is equivalent to \eqref{def:V}:
\begin{align*}
V = \bigg\{ \vvec \in \left( L^s(\Omega) \right)^n \Big| &
    (\Div \vvec)|_K \in L^{r'}(K) \quad \forall ~ K \in \Th \, , \\
& \displaystyle \sum_{K \in \Th} 
    \int_{\partial K} (\vvec \cdot \nvec_K) \, \varphi \, d\sigma = 0
    \quad \forall ~ \varphi \in \mathcal{D}(\Omega) \bigg\} \: .
\end{align*}
Here $\nvec_K$ denotes the exterior unit normal of $K$
and $\mathcal{D}(\Omega) = C_0^\infty(\Omega)$ is the space of test functions.
\end{lemma}
\begin{proof}
Evidently, ``$\subseteq$'' holds, since
$V \subset \left( L^s(\Omega) \right)^n$ and
$(\Div \vvec)|_K \in L^{r'}(K)$ for all $K \in \Th$.
Furthermore, for all $\vvec \in V$ and $\varphi \in \mathcal{D}(\Omega)$
Green's formula implies
\begin{align*}
\sum_{K \in \Th} \int_{\partial K} (\vvec \cdot \nvec_K) \, \varphi \, d\sigma
& = \sum_{K \in \Th} \left( \int_K \vvec \cdot \nabla \varphi \, d\xvec +
                              \int_K \Div \vvec \, \varphi \, d\xvec \right) \\
& = \int_\Omega \vvec \cdot \nabla \varphi \, d\xvec +
    \int_\Omega \Div\vvec \, \varphi \, d\xvec 
  = \int_{\partial \Omega} (\vvec \cdot \nvec) \, \varphi \, d\sigma = 0 \: .
\end{align*}
To prove that ``$\supseteq$'' holds, too, we must show that each
$\vvec \in \left( L^s(\Omega) \right)^n$, which fulfills the above requirements,
has a generalized divergence in $L^{r'}(\Omega)$.
To this end, we define $w \in L^{r'}(\Omega)$ by $w|_K = \Div \vvec$ for
all $K \in \Th$. Then for all $\varphi \in \mathcal{D}(\Omega)$ it holds
\begin{align*}
\int_\Omega w \varphi \, d\xvec & = 
\sum_{K \in \Th} \int_K \Div \vvec \, \varphi \, d\xvec = \sum_{K \in \Th} 
\left( \int_{\partial K} (\vvec \cdot \nvec_K) \, \varphi \, d\sigma -
       \int_K \vvec \cdot \nabla \varphi \, d\xvec \right) \\
& = - \sum_{K \in \Th} \int_K \vvec \cdot \nabla \varphi \, d\xvec =
- \int_\Omega \vvec \cdot \nabla \varphi \, d\xvec \: .
\end{align*}
Hence, $w$ is the (generalized) divergence of $\vvec$, as desired.
\end{proof}

Since every $\vvec_h|_K \in RT_0(K)$ is uniquely determined by the values
for the fluxes $v_e := \int_e \vvec_h|_K \cdot \nvec_e \, d\sigma$
across the edges $e \subset \partial K$ 
(cf.\ \cite[Prop. III.3.3]{Brezzi/Fortin}),
the fluxes $v_e$, $e \in \Eh$, can be used as degrees of freedom in
$V_h = RT_0(\Omega,\Th)$.
The corresponding basis functions $\left\{ \wvec_e \right\}_{e \in \Eh}$
are defined by the requirements
$$ \int_f \wvec_e \cdot \nvec_f \, d\sigma = \delta_{e,f} =
   \left\{ \begin{array}{ccc} 
		1 \: , & \text{if} & e=f \\
                0 \: , & \text{if} & e \neq f 
	   \end{array} \right. . $$

%%%%%%%%%%%%%%%%%%%%%%%%%%%%%%%%%%%%%%%%%%%%%%%%%%%%%%%%%%%%%%%%%%%%%%
%
\subsection{Hybridization}
\label{subsec:hybr}
The discrete mixed formulation reads:\\
Find $(\uvec_h,p_h) \in V_h \times Q_h$, such that
\begin{subequations}			\label{disc.var.pr}
\begin{align}				\label{disc.var.pr:a}
a(\uvec_h,\vvec_h) - b(\vvec_h,p_h) & = g(\vvec_h) \quad 
\text{for all } \vvec_h \in V_h \: ,\\
c(p_h,q_h) + b(\uvec_h,q_h) & = f(q_h) \quad \: 
\text{for all } q_h \in Q_h \: . 	\label{disc.var.pr:b}
\end{align}
\end{subequations}
For linear problems, the system of linear equations resulting 
from (\ref{disc.var.pr}) is difficult to solve.
Therefore this system is transformed by means of hybridization.
To this end, the finite-dimensional space $V_h$ is enlarged to the space $W_h$,
which is defined as follows:
$$ W_h := RT_{-1}(\Omega,\Th) =
   \big\{ \vvec_h \in \left( L^s(\Omega) \right)^n \bigm| 
           \vvec_h |_K \in RT_0(K) \text{ for all } K \in \Th \big\} \: . $$
Weakening the requirement $\vvec_h \in V$ to
$\vvec_h \in \left( L^s(\Omega) \right)^n$ means
that the normal components of the fluxes across interior edges must not
be continuous anymore.
Therefore $\vvec_h \in W_h$ is uniquely defined by the degrees of freedom
$ v_{K,e} := \int_e \vvec_h|_K \cdot \nvec_K \, d\sigma$
for $K \in \Th$ and $e \subset \partial K$.
Here $\nvec_K$ is the exterior unit normal with respect to $K$.
We denote the corresponding basis functions by $\wvec_{K,e}$.

To assure continuity of the fluxes, an additional equation is introduced,
where Lagrange-multipliers $\lambda_h \in \Lambda_h$ are employed.
Here $\Lambda_h$ is defined by
$$ \Lambda_h := 
   \left\{ \lambda_h \in L^2 ( E_h ) \bigm|
	   \lambda_h|_e \in P_0(e) \quad \text{for all } e \in \Eh \right\} , $$
where $E_h := \bigcup \{e \,|\, e \in \Eh\}$.
In particular, we consider subspaces $\Lambda_h^g$ of $\Lambda_h$, 
where for arbitrary $g \in L^1 (\partial\Omega)$
$$ \Lambda_h^{g} := 
   \left\{ \lambda_h \in \Lambda_h \Bigm|
	   \int_e (\lambda_h-g) \, d\sigma = 0
           \quad \text{for all } e \in \Eh^b \right\} \: . $$

Then the hybridized mixed formulation reads:\\
Find $(\uvec_h,p_h,\mu_h) \in W_h \times Q_h \times \Lambda_h^{g}$, 
such that
\begin{subequations}			\label{hybr.mixed}
\begin{align}				\label{hybr.mixed:a}
a(\uvec_h,\vvec_h) - \sum_{K \in \Th} b_K(\vvec_h,p_h) 
+ \sum_{K \in \Th} d_K(\mu_h,\vvec_h) & = 0
\hspace{2.5em} \text{for all } \vvec_h \in W_h \: ,\\
c(p_h,q_h) + \sum_{K \in \Th} b_K(\uvec_h,q_h) & = f(q_h)  
\hspace{0.7em} \text{for all } q_h \in Q_h \: , 	\label{hybr.mixed:b} \\
\sum_{K \in \Th} d_K(\lambda_h,\uvec_h) & = 0
\hspace{2.5em} \text{for all } \lambda_h \in \Lambda_h^0 \: ,	\label{hybr.mixed:c}
\end{align}
\end{subequations}
where $b_K(\vvec_h,q_h) := \int_K \Div \vvec_h q_h \, d\xvec$
and $d_K(\lambda_h,\vvec_h) := 
     \int_{\partial K} \lambda_h (\vvec_h \cdot \nvec_K) \, d\sigma$
for $K \in \Th$.
Note that we are allowed to identify the solutions $\uvec_h \in W_h$ and 
$p_h \in Q_h$ of (\ref{hybr.mixed}) with the solutions of (\ref{disc.var.pr}).

\begin{theorem}				\label{satz:hybr}
\begin{enumerate}	\renewcommand{\labelenumi}{\alph{enumi})}
\item Let $(\uvec_h,p_h) \in V_h \times Q_h$ be a solution of 
      \textup{(\ref{disc.var.pr})}. Then there exists a unique
      $\mu_h \in \Lambda_h^{g}$ such that
      $(\uvec_h,p_h,\mu_h)$ is a solution of \textup{(\ref{hybr.mixed})}.
\item Let $(\uvec_h,p_h,\mu_h) \in W_h \times Q_h \times \Lambda_h^{g}$
      be a solution of \textup{(\ref{hybr.mixed})}. 
      Then $(\uvec_h,p_h)$ is a solution of \textup{(\ref{disc.var.pr})}.
\end{enumerate}
\end{theorem}
Note that the nonlinear forms $a$ and $c$ remain unchanged in both formulations
(\ref{disc.var.pr}) and (\ref{hybr.mixed}).
Thus the proof of Theorem \ref{satz:hybr} does not differ from the linear case 
(cf.\ the derivation of Thm.~V.1.1 in \cite{Brezzi/Fortin}).

%%%%%%%%%%%%%%%%%%%%%%%%%%%%%%%%%%%%%%%%%%%%%%%%%%%%%%%%%%%%%%%%%%%%%%
%
\section{Equivalence to a nonconforming finite element method}
\label{sec:equiv}
%
%%%%%%%%%%%%%%%%%%%%%%%%%%%%%%%%%%%%%%%%%%%%%%%%%%%%%%%%%%%%%%%%%%%%%%
Like in \cite{Arbogast/Chen:95}, we can develop a nonconforming finite element
method for the numerical solution of equation (\ref{ell.eq}),
which yields the same solutions as (\ref{hybr.mixed}).
We consider the nonconforming ansatz space 
$\Nh := CR_1^0(\Th) \oplus \Bh$, where
\begin{align*}
CR_1^0(\Th) := \bigg\{ \xi_h \in L^2(\Omega) \, \Big| & \,
                       \xi_h|_K \in P_1(K) ~ \forall \: K \in \Th \, , \\[-1.5ex]
& \, \xi_h \text{ is continuous in } M_e ~ \forall \: e \in \Eh^i \, , \:
  \int_e \xi_h \, d\sigma = 0 ~ \forall \: e \in \Eh^b \bigg\}
\end{align*}
is the nonconforming Crouzeix--Raviart ansatz space of piecewise linear 
functions, which are continuous only at the midpoints $M_e$ of interior edges.
It is well known that $CR_1^0(\Th) \not\subset W^{1,1}(\Omega)$.
Since the midpoint rule is exact for affine-linear functions, i.e.\ 
$ \int_e \xi_h \, d\sigma = |e| \, \xi_h(M_e) $,
the value of the integral $\int_e \xi_h \, d\sigma$ is uniquely defined
for every interior edge $e \in \Eh^i$.
Moreover, we define the space of bubble functions
$$ \Bh := \left\{ \xi_h \in L^2(\Omega) \bigm|
   \xi_h|_K \in P_3(K) ~ \forall \: K \in \Th \, , \:
   \xi_h|_{\partial K} = 0 ~ \forall \: K \in \Th \right\} \: . $$
A basis for $\Bh$ is provided by
$\big\{ \prod_{i=1}^3 \lambda_i^K \bigm| K \in \Th \big\}$,
where $\lambda_i^K$ are the barycentric coordinates of the triangle $K$,
which are extended by zero outside of $K$.

To describe the nonconforming finite element method,
we need some projection operators.
At first we introduce the well-known $L^2$-projection operators 
$P_{Q_h}$, $P_{\Lambda_h}$ and $P_{\Lambda_h^0}$, which satisfy
$$ \begin{array}{c@{\,:\;}l@{~:\quad}rl}
P_{Q_h} & L^2(\Omega) \to Q_h & \displaystyle
\int_\Omega \left( \varphi - P_{Q_h} \varphi \right) q_h \, d\xvec = 0 &
\quad \text{for all } q_h \in Q_h \: , \\[2ex]
P_{\Lambda_h} & L^2(E_h) \to \Lambda_h & \displaystyle
\int_{E_h} \left( \varphi - P_{\Lambda_h} \varphi \right) 
                  \lambda_h \, d\sigma = 0 &
\quad \text{for all } \lambda_h \in \Lambda_h \: , \\[2ex]
P_{\Lambda_h^0} & L^2(E_h^i) \to \Lambda_h^0 & \displaystyle
\int_{E_h^i} \left( \varphi - P_{\Lambda_h^0} \varphi \right) 
                  \lambda_h \, d\sigma = 0 &
\quad \text{for all } \lambda_h \in \Lambda_h^0 \: .
\end{array} $$
Here $E_h^i := \bigcup \{e \,|\, e \in \Eh^i\}$.
Remember that the definitions of the spaces $Q_h$, $\Lambda_h$ 
and $\Lambda_h^0$ contain no continuity requirements.
Thus $P_{Q_h}$ is defined element by element,
$P_{\Lambda_h}$ and $P_{\Lambda_h^0}$ are defined edge by edge.

Using these projection operators, we can define a mapping
$\Sh : \Nh \to Q_h \times \Lambda_h^0$ by means of
$$ \Sh(\xi_h) = \left( P_{Q_h} \xi_h , P_{\Lambda_h^0} \xi_h \right) 
   \quad \text{for } \xi_h \in \Nh \: . $$
Since $\ker \Sh = 0$, $\Sh$ is an isomorphism (cf.\ \cite[(1.15)]{Brenner:92}).
In the following, let $\tilde{g}$ be a sufficiently smooth function on $\Omega$,
which satisfies $\int_e \tilde{g} \, d\sigma = \int_e g \, d\sigma$
for all $e \in \Eh^b$.
Then we can define a bijective mapping
$$ \widetilde{\mathcal{S}}_h^{\tilde{g}} : 
   \Nh+\tilde{g} \to Q_h \times \Lambda_h^{g} \quad , \quad
   \psi_h \mapsto \left( P_{Q_h} \psi_h , P_{\Lambda_h} \psi_h \right) \: . $$
Here, sufficiently smooth means that $P_{\Lambda_h} \psi_h$ has to be 
well defined for $\psi_h \in \Nh+\tilde{g}$.

Finally, we need a nonlinear projector
$\widetilde{P}_{W_h}: \left( L^s(\Omega) \right)^2 \to W_h$.
Generalizing the linear case considered in \cite{Arbogast/Chen:95}, we define
$\widetilde{P}_{W_h} \uvec$ for $\uvec \in \left( L^s(\Omega) \right)^2$
by the requirement
$$ a \big( \widetilde{P}_{W_h} \uvec , \vvec_h \big) 
  = \int_\Omega G \big( \widetilde{P}_{W_h}\uvec \big) \cdot \vvec_h \, d\xvec
  = a (\uvec , \vvec)
  \quad \text{for all } \vvec_h \in W_h \: . $$
Since $G$ is continuous, coercive and strictly monotone,
the Theorem of Browder and Minty \cite[Thm.~26.A]{Zeidler}
implies that $\widetilde{P}_{W_h} \uvec \in W_h$ is uniquely defined.
Like $P_{Q_h}$ also $\widetilde{P}_{W_h}$ is defined element by element,
since the definition of $W_h$ contains no continuity requirements.

Now, we can formulate the nonconforming finite element method:\\
Find $\psi_h \in \Nh + \tilde{g}$, such that for all $\xi_h \in \Nh$
\begin{equation}				\label{nonconf.FEM}
c(P_{Q_h} \psi_h, P_{Q_h} \xi_h) + \sum_{K \in \Th} 
  \int_K \widetilde{P}_{W_h} G^{-1}(\nabla \psi_h) \cdot \nabla \xi_h \, d\xvec =
f (P_{Q_h} \xi_h) \: .
\end{equation}
The following two lemmas, which establish the equivalence between 
(\ref{hybr.mixed}) and (\ref{nonconf.FEM}), generalize Thm.~1 
in \cite{Arbogast/Chen:95} to the nonlinear case.

\begin{lemma}				\label{lemma:hybrid=>nonconf}
Assume that $G(-\vvec) = -G(\vvec)$.
Let $\left( \uvec_h, p_h, \mu_h \right) \in W_h \times Q_h \times \Lambda_h^g$
be a solution of \textup{(\ref{hybr.mixed})}
and define $\psi_h \in \Nh + \tilde{g}$ by
$\psi_h := \big( \widetilde{\mathcal{S}}_h^{\tilde{g}} \big)^{-1} 
	   \big( (p_h,\mu_h) \big)\,$.
Then $\psi_h$ is a solution of \textup{(\ref{nonconf.FEM})} and
\begin{equation}					\label{u_h=}
\uvec_h = - \widetilde{P}_{W_h} G^{-1}(\nabla \psi_h) \: .
\end{equation}
\end{lemma}
\begin{proof}
We begin with the proof of \eqref{u_h=}. Let $\vvec_h \in W_h$ be given.
Since $(\Div \vvec_h)|_K \in P_0(K)$ and
$(\vvec_h|_K \cdot \nvec_K) |_e \in P_0(e)$ for every $e \subset \partial K$,
we obtain
%(see \cite[Prop.~III.3.2]{Brezzi/Fortin})
\begin{align*}
a \big( \widetilde{P}_{W_h} G^{-1}(\nabla \psi_h) , \vvec_h \big)
& = \sum_{K \in \Th} \int_K G \big( G^{-1}(\nabla \psi_h) \big)
                    \cdot \vvec_h \, d\xvec
 = \sum_{K \in \Th} \int_K \nabla \psi_h \cdot \vvec_h \, d\xvec \\
& = \sum_{K \in \Th} 
    \left( \int_{\partial K} \psi_h (\vvec_h \cdot \nvec_K) \, d\sigma
          - \int_K \Div \vvec_h \, \psi_h \, d\xvec \right) \\
& = \sum_{K \in \Th} 
    \left( \int_{\partial K} P_{\Lambda_h^0} \psi_h (\vvec_h \cdot \nvec_K) \, d\sigma
          - \int_K \Div \vvec_h \, P_{Q_h} \psi_h \, d\xvec \right) \\
& = \sum_{K \in \Th} 
    \left( \int_{\partial K} \mu_h (\vvec_h \cdot \nvec_K) \, d\sigma
          - \int_K \Div \vvec_h \, p_h \psi_h \, d\xvec \right) \\
& = \sum_{K \in \Th} d_K(\mu_h,\vvec_h) -
      \sum_{K \in \Th} b_h(\vvec_h,p_h)
  = - a(\uvec_h,\vvec_h) \\
& = - \int_\Omega G(\uvec) \cdot \vvec \, d\xvec
  = \int_\Omega G(-\uvec) \cdot \vvec \, d\xvec
  = a(-\uvec_h,\vvec_h) \: .
\end{align*}
Since $\widetilde{P}_{W_h} G^{-1}(\nabla \psi_h) \in W_h$ is uniquely defined,
\eqref{u_h=} follows.

To show that $\psi_h$ is a solution of \textup{(\ref{nonconf.FEM})},
we subtract (\ref{hybr.mixed:c}) from (\ref{hybr.mixed:b}).
Owing to the definitions of $b_K$ and $d_K$, this yields:
$$ c(p_h,q_h)
  + \sum_{K \in \Th} \int_K \! \Div \uvec_h \, q_h \, d\xvec
  - \sum_{K \in \Th} 
    \int_{\partial K} \! \lambda_h (\uvec_h \cdot \nvec_K) \, d\sigma
  = f(q_h) $$
for all $\left( q_h,\lambda_h \right) \in Q_h \times \Lambda_h^0$.
Next we define $\xi_h := \Sh^{-1} \left( (q_h,\lambda_h) \right)$
such that $q_h = P_{Q_h} \xi_h$ and $\lambda_h = P_{\Lambda_h^0} \xi_h$.
Hence we obtain
\begin{equation}				\label{proof:nonconf.FEM}
\begin{array}{r@{\;}c@{\;}l}
\displaystyle c(P_{Q_h} \psi_h , P_{Q_h} \xi_h) 
\displaystyle + \sum_{K \in \Th} \int_K \Div \uvec_h \, P_{Q_h} \xi_h \, d\xvec
& & \\[3ex]
\displaystyle - \sum_{K \in \Th} \int_{\partial K} P_{\Lambda_h^0} \xi_h \,
			(\uvec_h \cdot \nvec_K) \, d\sigma 
& = & f (P_{Q_h} \xi_h)
\end{array}
\end{equation}
for all $\xi_h \in \Nh$.
Employing (\ref{u_h=}) and the definitions of $P_{Q_h}$ and $P_{\Lambda_h^0}$,
the second and third term on the left hand side of (\ref{proof:nonconf.FEM}),
can be transformed as follows:
\begin{align*}
& \sum_{K \in \Th} \int_K \Div \uvec_h \, P_{Q_h} \xi_h \, d\xvec
- \sum_{K \in \Th} \int_{\partial K} P_{\Lambda_h^0} \xi_h
                                     (\uvec_h \cdot \nvec_K) \, d\sigma \\
& \hspace{9em} = \sum_{K \in \Th} \int_K \Div \uvec_h \, \xi_h \, d\xvec
   - \sum_{K \in \Th} \int_{\partial K} \xi_h \, 
                                       (\uvec_h \cdot \nvec_K) \, d\sigma \\
& \hspace{9em} = \sum_{K \in \Th} \int_K \!- \uvec_h \cdot \nabla \xi_h \, d\xvec
  = \sum_{K \in \Th} \int_K \! \widetilde{P}_{W_h} G^{-1}(\nabla \psi_h)
        		    \cdot \nabla \xi_h \, d\xvec \: .
\end{align*}
Inserting this identity into (\ref{proof:nonconf.FEM}) finally yields 
(\ref{nonconf.FEM}).
\end{proof}

\begin{lemma}				\label{lemma:nichtkonf=>hybrid}
Assume that $G(-\vvec) = -G(\vvec)$.
Let $\psi_h \in \Nh + \tilde{g}$ be a solution of \eqref{nonconf.FEM}.
Let $\uvec_h \in W_h$ be defined by \eqref{u_h=} and
$\left( p_h, \mu_h \right) \in Q_h \times \Lambda_h^g$ by
$\left( p_h, \mu_h \right) = \widetilde{\mathcal{S}}_h^{\tilde{g}} (\psi_h)$.
Then $\left( \uvec_h, p_h, \mu_h \right)$ is a solution of \eqref{hybr.mixed}.
\end{lemma}
\begin{proof}
Employing the definitions of $\uvec_h$, $p_h$ and $\mu_h$, we obtain 
\begin{align*}
a(\uvec_h,\vvec_h)
& = - a \big( \widetilde{P}_{W_h} G^{-1}(\nabla \psi_h) , \vvec_h \big)
 = - \sum_{K \in \Th} \int_K \nabla \psi_h \cdot \vvec_h \, d\xvec \\
& = \sum_{K \in \Th} \left( \int_K \Div \vvec_h \, \psi_h \, d\xvec 
    - \int_{\partial K} \psi_h (\vvec_h \cdot \nvec_K) \, d\sigma \right) \\
& = \sum_{K \in \Th} \int_K \Div \vvec_h \, P_{Q_h} \psi_h \, d\xvec 
 - \sum_{K \in \Th} \int_{\partial K} P_{\Lambda_h} \psi_h
				      (\vvec_h \cdot \nvec_K) \, d\sigma \\
& = \sum_{K \in \Th} \int_K \Div \vvec_h \, p_h \, d\xvec 
 - \sum_{K \in \Th} \int_{\partial K} \tilde{\mu}_h 
				      (\vvec_h \cdot \nvec_K) \, d\sigma \\
& = \sum_{K \in \Th} b_K(\vvec_h,p_h) - \sum_{K \in \Th} d_K(\mu_h,\vvec_h)
\end{align*}
for all $\vvec_h \in W_h$.
Hence (\ref{hybr.mixed:a}) is fulfilled.

Next we show that (\ref{hybr.mixed:b}) and (\ref{hybr.mixed:c}) are fulfilled, too.
Using (\ref{u_h=}) and partial integration, we obtain
\begin{align*}
\sum_{K \in \Th} \int_K \! \widetilde{P}_{W_h} G^{-1}(\nabla \psi_h) 
			\cdot \nabla \xi_h \, d\xvec 
& = \sum_{K \in \Th} \int_K \! - \uvec_h \cdot \nabla \xi_h \, d\xvec \\
& = \sum_{K \in \Th} \int_K \! \Div \uvec_h \, \xi_h \, d\xvec
  - \sum_{K \in \Th} \int_{\partial K} \! (\uvec_h \cdot \nvec_K) \xi_h \, d\sigma
\end{align*}
for all $\xi_h \in \Nh$.
Inserting this identity into \eqref{nonconf.FEM}
and using $p_h = P_{Q_h} \psi_h$ yields
$$ %\begin{align*}
c(p_h, P_{Q_h} \xi_h)
+ \sum_{K \in \Th} \int_K \Div \uvec_h \, P_{Q_h} \xi_h \, d\xvec %\hspace{10em}& \\
- \sum_{K \in \Th} \int_{\partial K} (\uvec_h \cdot \nvec_K) \,
				P_{\Lambda_h^0} \xi_h \, d\sigma
 = f(P_{Q_h} \xi_h) \: .
$$ %\end{align*}
Choosing $\xi_h$, such that $P_{Q_h} \xi_h = q_h \in Q_h$
and $P_{\Lambda_h^0} \xi_h = 0$, we obtain (\ref{hybr.mixed:b}).
In a similar manner, we obtain (\ref{hybr.mixed:c}),
if we choose $\xi_h$ with $P_{Q_h} \xi_h = 0$ and 
$P_{\Lambda_h^0} \xi_h = \lambda_h \in \Lambda_h^0$.
\end{proof}

%%%%%%%%%%%%%%%%%%%%%%%%%%%%%%%%%%%%%%%%%%%%%%%%%%%%%%%%%%%%%%%%%%%%%%
%
\section{Static Condensation}
\label{sec:stat.cond}
%
%%%%%%%%%%%%%%%%%%%%%%%%%%%%%%%%%%%%%%%%%%%%%%%%%%%%%%%%%%%%%%%%%%%%%%
Hybridization introduces many additional degrees of freedom into the problem.
Since most degrees of freedom can be eliminated locally,
this drawback can be avoided using static condensation.

The unknown functions $\uvec_h \in W_h$, $p_h \in Q_h$ and 
$\mu_h \in \Lambda_h^g$ can be represented by
$$ \uvec_h = \sum_{K \in \Th} \sum_{e \in \partial K} u_{K,e} \wvec_{K,e} \; ,
   \quad p_h = \sum_{K \in \Th} p_K \chi_K \; , \quad
   \mu_h = \sum_{e \in \Eh} \mu_e \chi_e \; . $$
Here $\wvec_{K,e}$ denote the basis functions of $W_h$
(cf.\ Subsection~\ref{subsec:hybr}), and $\chi_K$ and $\chi_e$ denote 
the characteristic functions of $K$ and $e$, resp.

Employing these basis functions in (\ref{hybr.mixed}) as test functions, i.e.,
$\vvec_h = \wvec_{K,\bar{e}}$ for $(K,\bar{e}) \in \Th \times \Eh$ 
with $\bar{e} \subset \partial K$,
$q_h = \chi_K$ for $K \in \Th$ and $\lambda_h = \chi_{\bar{e}}$ for 
$\bar{e} \in \Eh^i$, we obtain the following system of algebraic equations:
%for the degrees of freedom.
\begin{subequations}				\label{alg.Gls}
\begin{align}	
F_{K,\bar{e}} := & \;				\label{F_(K,e)}
a(\uvec_K , \wvec_{K,\bar{e}}) - p_K + \mu_{\bar{e}} = 0 \: , 
& K \in \Th , \: \bar{e} \subset \partial K \: , \\
F_K := & \;					\label{F_K}
c(p_K,\chi_K)
+ \sum_{e \subset \partial K} u_{K,e} - \int_K f \, d\xvec = 0 \: , 
\hspace{-5em} & K \in \Th \: , \\
F_{\bar{e}} := & \; u_{K,\bar{e}} + u_{K',\bar{e}} = 0 \: , 
& \bar{e} \in \Eh^i , \: \bar{e} = \partial K \cap \partial K' \: . \label{F_e}
\end{align}
\end{subequations}
Here $\uvec_K \in RT_0(K)$ is defined by 
$\uvec_K := \sum_{e \in \partial K} u_{K,e} \wvec_{K,e}$.
For the derivation of (\ref{alg.Gls}) we employed the fact 
that the basis functions $\wvec_{K,e}$ of $W_h$ satisfy
$$ \int_\Omega \Div \wvec_{K,e} \, d\xvec =
   \int_K \Div \wvec_{K,e} \, d\xvec =
   \int_{\partial K} \wvec_{K,e} \cdot \nvec_K \, d\sigma =
   \int_e \wvec_{K,e} \cdot \nvec_K \, d\sigma = 1 \: . $$

In this section we discuss three approaches to eliminate degrees of freedom
locally (i.e., on a single element).
This will lead to an enormous reduction of the number of equations
that have to be solved globally.

%%%%%%%%%%%%%%%%%%%%%%%%%%%%%%%%%%%%%%%%%%%%%%%%%%%%%%%%%%%%%%%%%%%%%%
%
\subsection{Elimination of the flux variables}
\label{subsec:Elim:m}
Considering \eqref{F_(K,e)}, we observe
that the degrees of freedom $u_{K,e}$ ($e \subset \partial K$) 
for fixed $K \in \Th$ appear only in equations (\ref{F_(K,e)})
for $\bar{e} \subset \partial K$.
Hence for each $K \in \Th$ we can eliminate $u_{K,e}$ 
for $e \subset \partial K$ by solving the system of equations
\begin{equation}					\label{impl.fct}
\Big( F_{K,\bar{e}} \big( \! \left( u_{K,e} \right)_{e \subset \partial K},
       p_K, \left( \mu_e \right)_{e \subset \partial K} \! \big) 
       \Big)_{\bar{e} \subset \partial K} = \vec{0} \: .
\end{equation}
In Lemma \ref{lemma:impl.fct} we show that this system has a unique solution.
For the remainder of this section we denote the triple of flux variables by 
$\vec{u}_K := \left( u_{K,e} \right)_{e \subset \partial K} \in \R^3$.
Since $\{ \wvec_{K,e} \, | \, e \subset \partial K \}$ 
is a basis of $V_K := RT_0(K)$, 
every $\vec{v}_K \in \R^3$ corresponds to a unique 
$\vvec_K = \sum\limits_{e \subset \partial K} v_{K,e} \wvec_{K,e} \in V_K$.
Thus we can identify $\vec{v}_K \in \R^3$ and $\vvec_K \in V_K$.

\begin{lemma}					\label{lemma:impl.fct}
Let $\big( p_K, (\mu_f)_{f \subset \partial K} \big) \in \R^4$ 
be given. Then there is a unique $\uvec_K \in V_K$,
that satisfies \textup{(\ref{impl.fct})}.
\end{lemma}

\begin{proof}
We equip $V_K = RT_0(K)$ with the $\left( L^s(K) \right)^2$-Norm
$\| \cdot \|_{V_K} = \| \cdot \|_{0,s,K}$.
System (\ref{impl.fct}) is equivalent to the variational problem
(cf.\ (\ref{hybr.mixed:a})):\\
Find $\uvec_K \in V_K$, such that
\begin{equation}				\label{hybr.mixed:K}
a_K(\uvec_K,\vvec_K) = 
\int_K \Div \vvec_K \, p_K \, d\xvec
- \sum_{f \subset \partial K} \int_{f} \mu_f (\vvec_K \cdot \nvec_K) \, d\sigma
=:\tilde{g}_K(\vvec_K) \: ,
\end{equation}
where $a_K$ denotes the nonlinear form defined by
$ a_K(\uvec_K,\vvec_K)
  = \int_K G(\uvec_K) \cdot \vvec_K \, d\xvec$.
Since $p_K$ and $\left( \mu_f \right)_{f \subset \partial K}$ are given,
$\tilde{g}_K$ is a continuous linear form on $V_K$.
We define an operator $A_K: V_K \to V_K'$ by means of
\begin{equation}					\label{A_K:def}
\left\langle A_K \uvec_K , \vvec_K \right\rangle_{V_K' \times V_K} 
:= a_K(\uvec_K,\vvec_K)
\quad \text{for } \uvec_K , \vvec_K \in V_K \: .
\end{equation}
Owing to the properties of $G$, $A_K$ is continuous, coercive 
and strictly monotone on $V_K$.
Therefore the Theorem of Browder and Minty \cite[Thm.~26.A]{Zeidler} yields
that there exists a unique solution $\uvec_K \in V_K$ of (\ref{hybr.mixed:K}).
\end{proof}

In general, the system (\ref{impl.fct}) has no closed form solution.
Then we can use Newton's method to find an approximative solution. 
To this end, it is necessary that the Jacobian matrix
$ \left( \partial F_{K,\bar{e}} / \partial u_{K,e} 
         \right)_{\bar{e},e \subset \partial K}$ 
is well-defined and invertible.

\begin{lemma}				\label{lemma:F_(K,e):diffable}
Let $G: \R^2 \to \R^2$ be continuously differentiable.
Then for all $K \in \Th$ and any choice of
$\big( p_K, ( \mu_f )_{f \subset \partial K} \big) \in \R^4$ 
the mapping
$$ \tilde{\mathbf{F}}_K : \R^3 \to \R^3 \quad , \quad
\left( u_{K,e} \right)_{e \subset \partial K} \mapsto
\left( F_{K,\bar{e}} \Big( (u_{K,e})_{e \subset \partial K}, 
			   p_K, ( \mu_f )_{f \subset \partial K} \Big)
       \right)_{\bar{e} \subset \partial K} $$
is continuously differentiable in $\R^3$.
%Here $F_{K,\bar{e}}$ is defined as in \eqref{F_(K,e)} and
%$\uvec_K = \sum_{e \in \partial K} u_{K,e} \wvec_{K,e}$.
\end{lemma}

\begin{proof}
We must show the continuity of the partial derivatives
$$ \frac{\partial F_{K,\bar{e}}}{\partial u_{K,\bar{f}}} (\vec{u}_K)
  = \frac{\partial f_{K,\bar{e}}}{\partial u_{K,\bar{f}}} (\vec{u}_K) \: ,
   \quad \bar{e},\bar{f} \subset \partial K \: , $$
where $f_{K,\bar{e}} : \R^3 \to \R$ is defined by
\begin{equation}				\label{f_(K,e)}
f_{K,\bar{e}} ( \vec{u}_K ) 
 := a_K(\uvec_K,\wvec_{K,\bar{e}})
  = \int_K G(\uvec_K) \cdot \wvec_{K,\bar{e}} \, d\xvec \: .
\end{equation}
To this end, it suffices to prove that the integrand
$G(\uvec_K) \cdot \wvec_{K,\bar{e}}$ is continuously differentiable.
This integrand is a composition of $G$ with the mappings
$\vec{u}_K \mapsto \uvec_K := \sum_{e \in \partial K} u_{K,e} \wvec_{K,e}$,
$\R^3 \to \R^2$, and $\xvec \mapsto \xvec \cdot \wvec_{K,\bar{e}}$, 
$\R^2 \to \R$.
Obviously, all these mappings are continuously differentiable.
Hence $\vec{u}_K \mapsto G(\uvec_K) \cdot \wvec_{K,\bar{e}}$, $\R^3 \to \R$,
is continuously differentiable, too.
\end{proof}

\begin{lemma}					\label{lemma:DF:invertible}
Assume that there is a constant $\underline{a} > 0$ 
and a norm $\| \cdot \|_{V_K}$ on $V_K$ such that 
the operator $A_K: V_K \to V_K'$ defined in \eqref{A_K:def} fulfills
%the inequality
\begin{equation}				\label{ass:DF:invertible}
\left\| A_K \uvec_K - A_K \vvec_K \right\|_{V_K'}
\ge \underline{a} \left\| \uvec_K - \vvec_K \right\|_{V_K} 
\quad \text{for all } \uvec_K, \vvec_K \in V_K \: .
\end{equation}
Then the Jacobian matrix
$$ D\tilde{\mathbf{F}}_K \left( \left( u_{K,e} \right)_{e \subset \partial K} \right) :=
   \left( \frac{\partial F_{K,\bar{e}}}{\partial u_{K,\bar{f}}}
          \left( \left( u_{K,e} \right)_{e \subset \partial K},
                 p_K, \left( \mu_f \right)_{f \subset \partial K} \right)   
   \right)_{\bar{e},\bar{f} \subset \partial K} $$
is invertible for all $\vec{u} = (u_{K,e})_{e \subset \partial K} \in \R^3$.
\end{lemma}

\begin{proof}
Since $V_K$ is finite-dimensional,
the euclidean norm $|\vec{v}_K|$ and $\|\vvec_K\|_{V_K}$ are equivalent.
In particular, there is a constant $C_1 > 0$ 
such that $|\vec{v}_K| \le C_1 \|\vvec_K\|_{V_K}$.
As observed in the proof of Lemma \ref{lemma:F_(K,e):diffable}, we have
$$ D\tilde{\mathbf{F}}_K (\vec{u}_K)
   = \left( \frac{\partial F_{K,\bar{e}}}{\partial u_{K,\bar{f}}} (\vec{u}_K)
          \right)_{\bar{e},\bar{f} \subset \partial K}
   = \left( \frac{\partial f_{K,\bar{e}}}{\partial u_{K,\bar{f}}} (\vec{u}_K)
            \right)_{\bar{e},\bar{f} \subset \partial K} \: , $$
where $f_{K,\bar{e}} : \R^3 \to \R$ is defined in (\ref{f_(K,e)}).
For every linear form $\vvec_K' \in V_K'$ there exists
a unique triple $\left( v'_{K,e} \right)_{e \subset \partial K} \in \R^3$
by means of evaluation at the test functions:
$$ \vvec_K' \mapsto \left( \left\langle \vvec_K', \wvec_{K,e} \right\rangle_{V_K' \times V_K} 
		    \right)_{e \subset \partial K}
   =: \left( v'_{K,e} \right)_{e \subset \partial K} $$ 
Hence $\|\vvec_K'\|_{V_K'}$ and 
$\big| \big( v'_{K,e} \big)_{e \subset \partial K} \big|$ 
are equivalent norms on $V_K'$.
In particular, there is a constant $C_2 > 0$ such that
$ \big\| \vvec_K' \big\|_{V_K'} \le 
  C_2 \, \big| \big( v'_{K,e} \big)_{e \subset \partial K} \big|\,$.

Now let the linear form $\vvec_K'$ be the image $A_K \uvec_K$ of
$\uvec_K = \sum\limits_{e \subset \partial K} u_{K,e} \wvec_{K,e} \in V_K$
for an arbitrary $\vec{u}_K \in \R^3$.
Then we have
$\left\langle A_K \uvec_K , \wvec_{K,\bar{e}} \right\rangle_{V_K' \times V_K}
 = f_{K,\bar{e}} \left( \vec{u}_K \right)$ for every 
$\bar{e} \subset \partial K$.
Therefore the following inequality holds for all $\vec{v}_K \in \R^3$:
\begin{equation}					\label{ineq}
\begin{array}{rcl}
\left| \vec{v}_K - \vec{u}_K \right| 
& \le & \displaystyle C_1 \left\| \vvec_K - \uvec_K \right\|_{V_K} 
  \le \frac{C_1}{\underline{a}} 
      \left\| A_K \vvec_K - A_K \uvec_K \right\|_{V_K'} \\[2ex]
& \le & \displaystyle \frac{C_1 C_2}{\underline{a}} 
        \left| \left( f_{K,\bar{e}} \left( \vec{v}_K \right) 
                     - f_{K,\bar{e}} \left( \vec{u}_K \right) 
               \right)_{\bar{e} \subset \partial K} \right| \: .
\end{array}
\end{equation}
Now it is easy to show that the Jacobian matrix 
$D\tilde{\mathbf{F}}_K \left( \vec{u}_K \right)$ is invertible.
Since 
$\left( f_{K,\bar{e}} \right)_{\bar{e} \subset \partial K}$ is differentiable,
we have
$$ \left( f_{K,\bar{e}} \left( \vec{v}_K \right) 
        - f_{K,\bar{e}} \left( \vec{u}_K \right) \right)_{
                                      \bar{e} \subset \partial K} 
  = D\tilde{\mathbf{F}}_K \left( \vec{u}_K \right) \left( \vec{v}_K - \vec{u}_K \right)
   + o \left( \left| \vec{v}_K - \vec{u}_K \right| \right) $$
for all $\vec{v}_K \in \R^3$. 
If $D\tilde{\mathbf{F}}_K \left( \vec{u}_K \right)$ was not invertible, 
there would exist $\vec{w}_K \in \R^3$ such that 
$D\tilde{\mathbf{F}}_K \left( \vec{u}_K \right) \vec{w}_K = \vec{0}$.
But this would imply that $\vec{v}_K = \vec{u}_K + \kappa \vec{w}_K$ 
for some $\kappa > 0$ would fulfill
$$ \left( f_{K,\bar{e}} \left( \vec{v}_K \right) 
        - f_{K,\bar{e}} \left( \vec{u}_K \right) \right)_{
					\bar{e} \subset \partial K} 
  = o \left( \left| \kappa \vec{w}_K \right| \right)
  = o \left( \left| \vec{v}_K - \vec{u}_K \right| \right) \: . $$
Consequently, we would obtain in the limit $\kappa \to 0$
$$ \lim_{|\vec{v}_K - \vec{u}_K| \to 0}
   \frac{ \left| \left( f_{K,\bar{e}} \left( \vec{v}_K \right) 
                       - f_{K,\bar{e}} \left( \vec{u}_K \right) 
                 \right)_{\bar{e} \subset \partial K} \right|}
        {\left| \vec{v}_K - \vec{u}_K  \right|} = 0 \: , $$
which is a contradiction to (\ref{ineq}).
\end{proof}

%%%%%%%%%%%%%%%%%%%%%%%%%%%%%%%%%%%%%%%%%%%%%%%%%%%%%%%%%%%%%%%%%%%%%%
%
\subsection{Elimination of flux and element variables -- 1st possibility}
\label{subsec:Elim:m&S:1}
Assuming additional monotonicity requirements for $R$, we can
%in addition to the elimination of the flux variables $u_{K,e}$ for 
%$e \subset \partial K$ by solving the systems of equations (\ref{impl.fct})
eliminate also the element variables $p_K$ by solving equations (\ref{F_K}).
A closer look onto (\ref{F_K}) shows that
$$ \int_K f \, d\xvec - \sum_{e \subset \partial K} u_{K,e}
   = c(p_K,\chi_K)
   = \int_\Omega R(p_K) \chi_K \, d\xvec 
   = \int_K R(p_K) \, d\xvec \: . $$
To ensure that $p_K$ is uniquely defined by this identity
we have to require that the mapping
$C_K: \R \to \R$, $p_K \mapsto \int_K R(p_K) \, d\xvec$ is strictly monotone.
%i.e., $(C_K p_K - C_K q_K)(p_k-q_K) > 0$ for $p_K \neq q_K$.
Then $C_K$ is invertible and we can solve the above equation:
$$ p_K = p_K ( \vec{u}_K ) :=
   C_K^{-1} \left( \int_K f \, d\xvec
	         - \sum_{e \subset \partial K} u_{K,e} \right) \, . $$
Inserting this relation into (\ref{F_(K,e)}) yields
\begin{equation}				\label{bar(F)_(K,e)}
\bar{F}_{K,\bar{e}}(\vec{u}_K)
 := a(\uvec_K,\wvec_{K,e}) - p_K ( \vec{u}_K ) + \mu_{\bar{e}}
  = f_{K,\bar{e}}(\vec{u}_K) - p_K ( \vec{u}_K ) + \mu_{\bar{e}} = 0 \: ,
\end{equation}
where $f_{K,\bar{e}}$ is still defined as in (\ref{f_(K,e)}).
Now we can solve the systems of equations
\begin{equation}				\label{bar(F)=0}
\bar{\mathbf{F}}_K (\vec{u}_K) :=
\left( \bar{F}_{K,\bar{e}}(\vec{u}_K) \right)_{\bar{e} \subset \partial K} = 
\vec{0}
\end{equation}
to eliminate the flux variables $u_{K,e}$ for $e \subset \partial K$.

\begin{lemma}					\label{lemma:bar(F):impl.fct}
Let $C_K$ be strictly monotone.
Then for arbitrary $( \mu_f )_{f \subset \partial K} \in \R^3$ 
there exists a unique triple 
$\vec{u}_K = \left( u_{K,e} \right)_{e \subset \partial K} \in \R^3$,
such that \eqref{bar(F)=0} is satisfied.
\end{lemma}

\begin{proof}
System (\ref{bar(F)=0}) is equivalent to the following variational problem 
(cf.\ the proof of Lemma \ref{lemma:impl.fct}):
Find $\uvec_K \in V_K$, such that for all $\vvec_K \in V_K$
\begin{equation}				\label{var.pr:bar(F)}
a_K(\uvec_K,\vvec_K)- b_K \big( \vvec_K , p_K(\vec{u}_K) \big)
= - d_K(\mu_K,\vvec_K) \: .
\end{equation}
Here we equip $V_K = RT_0(K)$ with the $\left( L^s(K) \right)^2$-norm
$\| \cdot \|_{V_K} = \| \cdot \|_{0,s,K}$.
If $\mu_K := \sum_{e \subset \partial K} \mu_e \chi_e$ is given,
we have $d_K(\mu_K,\cdot) \in V_K'$.
Let the operator $A_K$ be defined as in \eqref{A_K:def}.
Furthermore, we employ the (linear) operator $B_K: V_K \to Q_K'$, defined by 
$\langle B_K \vvec_K , q_K \rangle_{Q_K' \times Q_K} = b_K(\vvec_K , q_K)$
for $\vvec_K \in V_K$ and $q_K \in Q_K$, where $Q_K := P_0(K)$ is equipped 
with the norm $\| \cdot \|_{Q_K} = \| \cdot \|_{0,r,K}$.
Identifying $\vec{u}_K \in \R^3$ with $\uvec_K \in V_K$ again,
we define a nonlinear operator $\bar{A}_K: V_K \to V_K'$ by
\begin{equation}					\label{bar(A):def}
\big\langle \bar{A}_K \uvec_K , \vvec_K \big\rangle_{V_K' \times V_K}
:= \langle A_K \uvec_K , \vvec_K \rangle_{V_K' \times V_K} \!
  - \langle B_K \vvec_K , p_K(\vec{u}_K) \rangle_{Q_K' \times Q_K} .
\end{equation}
Evidently, $\bar{A}_K$ is continuous.
Furthermore, $\bar{A}_K$ inherits the strict monotonicity and coercivity
of $A_K$, since for $\uvec_K, \vvec_K \in V_K$
\begin{align*}
\langle B_K (\uvec_K - \vvec_K) , 
p_K(\vec{u}_K) - p_K(\vec{v}_K) \rangle_{Q_K' \times Q_K} = 
b_K(\uvec_K - \vvec_K , p_K(\vec{u}_K) - p_K(\vec{v}_K) ) & = \\
\int_K \Div ( \uvec_K - \vvec_K ) 
       \left( p_K(\vec{u}_K) -p_K(\vec{v}_K) \right) \, d\xvec \le & \: 0 \: .
\end{align*}
This inequality follows in the following manner from the monotonicity of $C_K$:
\begin{align*}
& \Div ( \uvec_K - \vvec_K ) 
	 \big( p_K(\vec{u}_K) - p_K(\vec{v}_K) \big)
= \left( \sum_e u_{K,e} - \sum_e v_{K,e} \right) 
      \big( p_K(\vec{u}_K) - p_K(\vec{v}_K) \big) \\
& = \: \left( \bigg( \int_K f \, d\xvec
	           - \sum_{e \subset \partial K} v_{K,e} \bigg) 
            - \bigg( \int_K f \, d\xvec
	            - \sum_{e \subset \partial K} u_{K,e} \bigg) \right)
      \big( p_K(\vec{u}_K) - p_K(\vec{v}_K) \big) \\
& = \: \big( C_K \big(p_K(\vec{v}_K)\big)
	    - C_K \big(p_K(\vec{u}_K)\big) \big)
      \big( p_K(\vec{u}_K) - p_K(\vec{v}_K) \big) \le 0 \: .
\end{align*}
Again, the Theorem of Browder and Minty \cite[Thm.~26.A]{Zeidler} yields
that there exists a unique solution $\uvec_K \in V_K$ of (\ref{var.pr:bar(F)}).
The corresponding triple $\vec{u} \in \R^3$ is the solution of (\ref{bar(F)=0}).
\end{proof}

Again, (\ref{bar(F)=0}) can not be solved in closed form.
Hence, we have to use a numerical method like Newton's method.
Under certain assumptions on $G$ and $R$, we can show that
$\bar{\mathbf{F}}_K$ is continuously differentiable
and that the Jacobian matrix $D \bar{\mathbf{F}}_K$ is invertible.

\begin{lemma}				\label{lemma:bar(F)_(K,e)diffable}
Let $G: \R^2 \to \R^2$ and $C_K:\R \to \R$ be continuously differentiable
and let $C_K' \neq 0$ on $\R$.
Then for all $K \in \Th$ and any choice of
$( \mu_f )_{f \subset \partial K} \in \R^3$ 
the mapping
$$ \bar{\mathbf{F}}_K : \R^3 \to \R^3 \quad , \quad
   \vec{u} = \left( u_{K,e} \right)_{e \subset \partial K} \mapsto
   \left( \bar{F}_{K,\bar{e}} 
          \left( \left( u_{K,e} \right)_{e \subset \partial K},
                \left( \mu_f \right)_{f \subset \partial K} \right) 
          \right)_{\bar{e} \subset \partial K} $$
is continuously differentiable in $\R^3$.
%Here $\bar{F}_{K,\bar{e}}$ is defined as in \eqref{bar(F)_(K,e)} and
%$\uvec_K = \sum_{e \in \partial K} u_{K,e} \wvec_{K,e}$.
\end{lemma}
\begin{proof}
Like in the proof of Lemma \ref{lemma:F_(K,e):diffable} we show that
the partial derivatives
$$ \frac{\partial \bar{F}_{K,\bar{e}}}{\partial u_{K,\bar{f}}} (\vec{u})
 = \frac{\partial f_{K,\bar{e}}}{\partial u_{K,\bar{f}}} (\vec{u})
  - \frac{\partial p_K}{\partial u_{K,\bar{f}}} (\vec{u}) $$
for $\bar{e},\bar{f} \subset \partial K$ exist and are continuous.
%Here $f_{K,\bar{e}}$ is defined in \eqref{f_(K,e)} and 
The partial derivatives $\partial f_{K,\bar{e}}/\partial u_{K,\bar{f}}$
are considered in the proof of Lemma \ref{lemma:F_(K,e):diffable}.
Thus we only have to deal with the partial derivatives
$\partial p_K/\partial u_{K,\bar{f}}$.
Since $p_K$ is the solution of (\ref{F_K}),
we can compute these derivatives by means of the implicit function theorem:
$$ \frac{\partial p_K}{\partial u_{K,\bar{f}}} (\vec{u})
   = - \left( \frac{\partial F_K}{\partial p_K} 
	      \left( \vec{u}, p_K(\vec{u}) \right) \right)^{-1} 
     \frac{\partial F_K}{\partial u_{K,\bar{f}}} 
	      \left( \vec{u}, p_K(\vec{u}) \right) \: . $$
Evidently, $\partial F_K/\partial u_{K,\bar{f}} = 1$ for all 
$\bar{f} \subset \partial K$ and $\partial F_K/\partial p_K = C_K'$.
Since $C_K' \neq 0$, $p_K(\vec{u})$ is continuously differentiable.
\end{proof}

\begin{lemma}				\label{lemma:bar(F)_(K,e):invertible}
Assume that the assumptions of Lemmas \textup{\ref{lemma:DF:invertible}},
\textup{\ref{lemma:bar(F):impl.fct}} 
and \textup{\ref{lemma:bar(F)_(K,e)diffable}} are fulfilled.
Then for all $K \in \Th$ and any choice of
$( \mu_f )_{f \subset \partial K} \in \R^3$ 
the Jacobian matrix
$$ D\bar{\mathbf{F}}_K(\vec{u}) :=
   \left( \frac{\partial \bar{F}_{K,\bar{e}}}{\partial u_{K,\bar{f}}}
          \left( \left( u_{K,e} \right)_{e \subset \partial K},
                 \left( \mu_f \right)_{f \subset \partial K} \right)   
   \right)_{\bar{e},\bar{f} \subset \partial K} $$
is invertible.
\end{lemma}
\begin{proof}
Like in the proof of Lemma \ref{lemma:bar(F):impl.fct} we can conclude that
\begin{align*}
\big\langle \bar{A}_K \uvec_K - \bar{A}_K \vvec_K , 
    \uvec_K - \vvec_K \big\rangle_{V_K' \times V_K} \: & \ge 
\langle A_K \uvec_K - A_K \vvec_K , \uvec_K - \vvec_K \rangle_{V_K' \times V_K} \\
& \ge \underline{a} \left\| \uvec_K - \vvec_K \right\|_{V_K}^2 
    \quad \text{for all } ~ \uvec_K , \vvec_K \in V_K \: .
\end{align*}
Hence we obtain in the same way as in the proof of Lemma \ref{lemma:DF:invertible} 
\begin{align*}
\left| \vec{v}_K - \vec{u}_K \right| \:
\le & \: C_1 \left\| \vvec_K - \uvec_K \right\|_{V_K} 
\le \frac{C_1}{\underline{a}} 
      \left\| \bar{A}_K \vvec_K - \bar{A}_K \uvec_K \right\|_{V_K'} \\[2ex]
\le & \: \frac{C_1 C_2}{\underline{a}} 
        \left| \left( (f_{K,\bar{e}} - p_K) ( \vec{v}_K ) 
                     - (f_{K,\bar{e}} - p_K) ( \vec{u}_K ) 
               \right)_{\bar{e} \subset \partial K} \right|
\end{align*}
for all $\vec{v}, \vec{u} \in \R^3$.
This implies the invertibility of $D\bar{\mathbf{F}}_K(\vec{u})$.
\end{proof}

%%%%%%%%%%%%%%%%%%%%%%%%%%%%%%%%%%%%%%%%%%%%%%%%%%%%%%%%%%%%%%%%%%%%%%
%
\subsection{Elimination of flux and element variables -- 2nd possibility}
\label{subsec:Elim:m&S:2}
There is another possibility to eliminate the flux variables $u_{K,e}$ for
$e \subset \partial K$ and the element variables $p_K$ for every $K \in \Th$.
We can consider the system %consisting of the four equations
\begin{equation}				\label{F_K=0}
\begin{array}{rcl}
F_{K,\bar{e}} \left( \left( u_{K,e} \right)_{e \subset \partial K}, p_K,
		     \left( \mu_e \right)_{e \subset \partial K} \right)
& = & 0 \; , \quad \bar{e} \subset \partial K \: , \\[1ex]
F_K \left( \left( u_{K,e} \right)_{e \subset \partial K}, p_K \right) & = & 0
\end{array}
\end{equation}
which defines (implicitly) the local variables $u_{K,e}$ 
($e \subset \partial K$) and $p_K$ depending on 
$\left( \mu_e \right)_{e \subset \partial K}$.
To this end we define for given values of $\mu_f$ ($f \subset \partial K$)
a mapping $\mathbf{F}_K : \R^4 \to \R^4$ via
\begin{equation}				\label{bf(F)_K}
\mathbf{F}_K \left( (u_{K,e})_{e \subset \partial K} , p_K \right) =
   \left( \!\! \begin{array}{c}
		\left( F_{K,\bar{e}} 
                   \left( \left( u_{K,e} \right)_{e \subset \partial K}, p_K ,
                          \left( \mu_f \right)_{f \subset \partial K} \right) 
                \right)_{\bar{e} \subset \partial K} \\
	    	F_K \left( \left( u_{K,e} \right)_{e \subset \partial K},
                           \left( \mu_f \right)_{f \subset \partial K} \right)
	  \end{array} \!\!
   \right) \, .
\end{equation}
Of course, (\ref{F_K=0}) is uniquely solvable, too:

\begin{lemma}
Let $R$ be strictly monotone. Then for any choice of
$( \mu_f )_{f \subset \partial K} \in \R^3$ there exists a unique solution
$\left( ( u_{K,e} )_{e \subset \partial K} , p_K \right) \in \R^4$ 
of \eqref{F_K=0}.
\end{lemma}
\begin{proof}
Evidently, system (\ref{F_K=0}) is equivalent to the variational problem:
Find $(\uvec_K,p_K) \in V_K \times Q_K$ such that 
\begin{align*}
a_K(\uvec_K,\vvec_K) \, d\xvec - b_K(\vvec_K,p_K) \: & = - d_K(\mu_K,\vvec_K)
& \quad \text{for all } \vvec_K \in V_K \: , \\
c_K(p_K,q_K) + b_K(\uvec_K,q_K) \: & = f(q_K)
& \quad \text{for all } q_K \in Q_K \: .
\end{align*}
The unique solvability of this problem follows similarly to 
\cite[Thm.~1.8 or Thm.~2.4]{Knabner/Summ}.
\end{proof}

To solve (\ref{F_K=0}) by Newton's method, we need again:

\begin{lemma}				\label{lemma:bf(F)_K:diffable}
Let the assumptions of Lemma \textup{\ref{lemma:bar(F)_(K,e)diffable}}
be fulfilled.
Then for every $K \in \Th$ and any choice of
$( \mu_f )_{f \subset \partial K} \in \R^3$ the mapping $\mathbf{F}_K$ 
from \eqref{bf(F)_K} is continuously differentiable in $\R^4$.
\end{lemma}
\begin{proof}
We have seen already in the proof of Lemma \ref{lemma:bar(F)_(K,e)diffable}
that the partial derivatives
$$ \frac{\partial F_{K,\bar{e}}}{\partial u_{K,\bar{f}}} ~ , \quad
   \frac{\partial F_{K,\bar{e}}}{\partial p_K} = -1 ~ , \quad
   \frac{\partial F_K}{\partial u_{K,\bar{f}}} = 1 ~ \text{and} \quad
   \frac{\partial F_K}{\partial p_K} = C_K' $$
exist and are continuous.
\end{proof}

\begin{lemma}				\label{lemma:bf(F)_K:invertible}
Let the assumptions of Lemma \textup{\ref{lemma:bar(F)_(K,e):invertible}}
are fulfilled.
Then for all $K \in \Th$ and any choice of
$( \mu_f )_{f \subset \partial K} \in \R^3$ 
the Jacobian matrix
$$ D\mathbf{F}_K \left( (u_{K,e})_e, p_K \right) :=
   \frac{\partial \mathbf{F}_K \left( (u_{K,e})_e, p_K \right)}
	{\partial \! \left( (u_{K,e})_e,p_K \right)} $$
is invertible.
\end{lemma}
\begin{proof}
We show that the determinant of
$$ D\mathbf{F}_K =
   \left( \begin{array}{cccc}
	  \frac{\partial F_{K,e_1}}{\partial u_{K,e_1}} &
	  \frac{\partial F_{K,e_1}}{\partial u_{K,e_2}} &
	  \frac{\partial F_{K,e_1}}{\partial u_{K,e_3}} &
	  -1 \\[1ex]
	  \frac{\partial F_{K,e_2}}{\partial u_{K,e_1}} &
	  \frac{\partial F_{K,e_2}}{\partial u_{K,e_2}} &
	  \frac{\partial F_{K,e_2}}{\partial u_{K,e_3}} &
	  -1 \\[1ex]
	  \frac{\partial F_{K,e_3}}{\partial u_{K,e_1}} &
	  \frac{\partial F_{K,e_3}}{\partial u_{K,e_2}} &
	  \frac{\partial F_{K,e_3}}{\partial u_{K,e_3}} &
	  -1 \\[1ex]
	  1 & 1 & 1 & \frac{\partial F_K}{\partial p_K}
	  \end{array} \right)
$$
does not vanish.
Since $\frac{\partial F_K}{\partial p_K} = C_K' > 0$, 
we can multiply the last row
by $\left(\frac{\partial F_K}{\partial p_K}\right)^{-1}$ 
and add the result successively to the first three rows.
This way, we eliminate the first three entries of the last column.
An expansion of the determinant along the last column finally yields
$$ \det(D\mathbf{F}_K) = 
   \frac{\partial F_K}{\partial p_K} \det(D\bar{\mathbf{F}}_K) \neq 0 \: , $$
because $D\bar{\mathbf{F}}_K$ is invertible,
as we have seen in the proof of Lemma \ref{lemma:bar(F)_(K,e)diffable}. 
\end{proof}

%%%%%%%%%%%%%%%%%%%%%%%%%%%%%%%%%%%%%%%%%%%%%%%%%%%%%%%%%%%%%%%%%%%%%%
%
\section{Solution of resulting system of equations}
\label{sec:sol}
%
%%%%%%%%%%%%%%%%%%%%%%%%%%%%%%%%%%%%%%%%%%%%%%%%%%%%%%%%%%%%%%%%%%%%%%
All the possibilities for static condensation presented in the last section
reduce the number of equations and of degrees of freedom.
If we eliminate only the flux variables, as described in Subsection 
\ref{subsec:Elim:m}, we still have to compute the values
$(p_K)_{K \in \Th}$ of element variables and
$(\mu_e)_{e \in \Eh^i}$ of Lagrange-multipliers. 
If we eliminate the flux and element variables, as described in Subsections
\ref{subsec:Elim:m&S:1} and \ref{subsec:Elim:m&S:2},
we still have to compute the values $(\mu_e)_{e \in \Eh^i}$ 
of Lagrange-multipliers.

%%%%%%%%%%%%%%%%%%%%%%%%%%%%%%%%%%%%%%%%%%%%%%%%%%%%%%%%%%%%%%%%%%%%%%
%
\subsection{Linearization by Newton's method}
\label{subsec:Newton}
The elimination of variables, as described in Section \ref{sec:stat.cond},
introduces nonlinearities into the originally linear equation (\ref{F_e}).
If we proceed as described in Subsection \ref{subsec:Elim:m},
we introduce also additional nonlinearities into (\ref{F_K}).
We propose Newton's method for the linearization of the resulting system 
of equations.
Below, we describe, how the entries of the Jacobian matrix can be computed.

%%%%%%%%%%%%%%%%%%%%%%%%%%%%%%%%%%%%%%%%%%%%%%%%%%%%%%%%%%%%%%%%%%%%%%%%%%%
\subsection*{Computation of the Jacobian matrix 
		after elimination of flux variables}
After solving the system of equations \eqref{impl.fct} for some $K \in \Th$, 
the flux variables $u_{K,e}$ ($e \subset \partial K$) are uniquely defined
as implicit functions depending on $p_K$ and $(\mu_f)_{f \in \partial K}$:
$$ u_{K,e} = u_{K,e} \left( p_K, (\mu_f)_{f \subset \partial K} \right) 
   \quad \text{for } e \in \partial K \: . $$
While computing the partial derivatives of $F_K$ and $F_{\bar{e}}$,
we have to take this dependence into account.

Hence we obtain the following expressions for the entries of the Jacobian 
matrix:
$$ \begin{array}{r@{\;}c@{\;}lr}
\displaystyle \frac{\partial F_K}{\partial p_K} 
& = & \lefteqn{\sum_{e \subset \partial K} 
      \frac{\partial u_{K,e}}{\partial p_K} 
      \left( p_K, \left( \mu_f \right)_{f \subset \partial K} \right) +
      C_K' \: ,}
    & \forall \; K \in \Th \, , \\[3ex]
\displaystyle \frac{\partial F_K}{\partial \mu_{\bar{f}}} 
& = & \lefteqn{\sum_{e \subset \partial K} 
      \frac{\partial u_{K,e}}{\partial \mu_{\bar{f}}} 
      \left( p_K, \left( \mu_f \right)_{f \subset \partial K} \right) \: ,}
    & \forall \; K \in \Th \, , \; \bar{f} \subset \partial K \, , \\[3ex]
\displaystyle \frac{\partial F_{\bar{e}}}{\partial p_K} 
& = & \displaystyle \frac{\partial u_{K,\bar{e}}}{\partial p_K} 
      \left( p_K, \left( \mu_f \right)_{f \subset \partial K} \right) \: ,
    & \forall \; \bar{e} \in \Eh^i \, , 
              \; \bar{e} \subset \partial K \, , \\[3ex]
\displaystyle \frac{\partial F_{\bar{e}}}{\partial \mu_{\bar{f}}}
& = & \lefteqn{\frac{\partial u_{K,\bar{e}}}{\partial \mu_{\bar{f}}}
      \left( p_K, \left( \mu_f \right)_{f \subset \partial K} \right)
      + \delta_{\bar{e},\bar{f}} 
      \frac{\partial u_{K',\bar{e}}}{\partial \mu_{\bar{f}}}
      \left( p_{K'}, \left( \mu_{f'} \right)_{f' \subset \partial K'} 
      \right) \: ,} \\[2ex]
& & & \qquad \forall \; \bar{e},\bar{f} \in \Eh^i \, ,
              \; \bar{e} = \partial K \cap \partial K' \, , 
              \; \bar{f} \subset \partial K \: .
\end{array} $$
Note that we need the partial derivatives 
$\partial u_{K,\bar{e}}/\partial p_K$ resp.\ 
$\partial u_{K,\bar{e}}/\partial \mu_{\bar{f}}$.
By means of the implicit function theorem these partial derivatives 
are given by:
\begin{align*}
\left( \frac{\partial u_{K,\bar{e}}}{\partial p_K} 
       \right)_{\bar{e} \subset \partial K} \:
& = - D\tilde{\mathbf{F}}_K^{-1}
\left( \frac{\partial F_{K,\bar{e}}}{\partial p_K} 
       \right)_{\bar{e} \subset \partial K} \\
\mbox{and } \hspace{20ex} & \hspace{58ex}\\
\left( \frac{\partial u_{K,\bar{e}}}{\partial \mu_{\bar{f}}} 
       \right)_{\bar{e},\bar{f} \subset \partial K} \:
& = - D\tilde{\mathbf{F}}_K^{-1}
\left( \frac{\partial F_{K,\bar{e}}}{\partial \mu_{\bar{f}}} 
       \right)_{\bar{e},\bar{f} \subset \partial K} \: . 
\end{align*}
We therefore need the inverse of the local matrix $D\tilde{\mathbf{F}}_K$.
In Lemma \ref{lemma:DF:invertible} we showed that this inverse does exist.
In addition, we need the following partial derivatives:
$$ \frac{\partial F_{K,\bar{e}}}{\partial p_K} = - 1 \quad \text{and} \quad
\frac{\partial F_{K,\bar{e}}}{\partial \mu_{\bar{f}}}
 = \delta_{\bar{e},\bar{f}} \: . $$

%%%%%%%%%%%%%%%%%%%%%%%%%%%%%%%%%%%%%%%%%%%%%%%%%%%%%%%%%%%%%%%%%%%%%%%%%%%
\subsection*{Computation of the Jacobian matrix 
		after elimination of flux and element variables 
		according to the 1st possibility}
After elimination of flux and element variables as described in Subsection
\ref{subsec:Elim:m&S:1}, we have to solve only the system
$\left( F_{\bar{e}} = 0 \right)_{\bar{e} \in \Eh^i}$ for the unknowns
$\left(\mu_e\right)_{e \in \Eh^i}$.
Since $u_{K,\bar{e}}$ depends nonlinearly on $\mu_e$ ($e \subset \partial K$),
these equations become nonlinear.
To solve this system by means of Newton's method 
we therefore need the partial derivatives
$$ \frac{\partial F_{\bar{e}}}{\partial \mu_{\bar{f}}}
   = \frac{\partial u_{K,\bar{e}}}{\partial \mu_{\bar{f}}}
      \left( \left( \mu_f \right)_{f \subset \partial K} \right)
    + \delta_{\bar{e},\bar{f}} 
     \frac{\partial u_{K',\bar{e}}}{\partial \mu_{\bar{f}}}
      \left( \left( \mu_{f'} \right)_{f' \subset \partial K'} \right) $$
for $\bar{e},\bar{f} \in \Eh^i$, $\bar{e} = \partial K \cap \partial K'$,
$\bar{f} \subset \partial K$.
Again, the implicit function theorem yields the formula for
the computation of the partial derivatives
$\partial u_{K,\bar{e}}/\partial \mu_{\bar{f}}$:
$$ \left( \frac{\partial u_{K,\bar{e}}}{\partial \mu_{\bar{f}}} 
          \right)_{\bar{e},\bar{f} \subset \partial K} 
  = - D\bar{\mathbf{F}}_K^{-1}
      \left( \frac{\partial \bar{F}_{K,\bar{e}}}{\partial \mu_{\bar{f}}} 
             \right)_{\bar{e},\bar{f} \subset \partial K} \: . $$
Since
$~ \displaystyle \frac{\partial \bar{F}_{K,\bar{e}}}{\partial \mu_{\bar{f}}}
          = \delta_{\bar{e},\bar{f}} \: , ~$
we have
$~ \displaystyle \left( \frac{\partial u_{K,\bar{e}}}{\partial \mu_{\bar{f}}} 
                 \right)_{\bar{e},\bar{f} \subset \partial K}
         = - D\bar{\mathbf{F}}_K^{-1}$.
The existence of $D\bar{\mathbf{F}}_K^{-1}$ is established 
in Lemma \ref{lemma:bar(F)_(K,e):invertible}.

%%%%%%%%%%%%%%%%%%%%%%%%%%%%%%%%%%%%%%%%%%%%%%%%%%%%%%%%%%%%%%%%%%%%%%%%%%%
\subsection*{Computation of the Jacobian matrix 
		after elimination of flux and element variables 
		according to the 2nd possibility}
After elimination of flux and element variables as described in Subsection
\ref{subsec:Elim:m&S:2}, we again have to solve the system
$\left( F_{\bar{e}} = 0 \right)_{\bar{e} \in \Eh^i}$ for
$\left(\mu_e\right)_{e \in \Eh^i}$.
Since these equations become nonlinear, we propose Newton's method
to solve them. Again, we need the partial derivatives
$$ \frac{\partial F_{\bar{e}}}{\partial \mu_{\bar{f}}}
   = \frac{\partial u_{K,\bar{e}}}{\partial \mu_{\bar{f}}}
      \left( \left( \mu_f \right)_{f \subset \partial K} \right)
    + \delta_{\bar{e},\bar{f}} 
     \frac{\partial u_{K',\bar{e}}}{\partial \mu_{\bar{f}}}
      \left( \left( \mu_{f'} \right)_{f' \subset \partial K'} \right) $$
for $\bar{e},\bar{f} \in \Eh^i$, $\bar{e} = \partial K \cap \partial K'$,
$\bar{f} \subset \partial K$.
Now the implicit function theorem yields the following formula for
the computation of these partial derivatives:
$$ \frac{\partial \left( ( u_{K,\bar{e}} )_{\bar{e} \subset \partial K}, 
                         p_K \right)}
        {\partial \left( ( \mu_{\bar{f}})_{\bar{f} \subset \partial K}\right)}
  = D\mathbf{F}_K^{-1}
    \frac{\partial \mathbf{F}_K}
	 {\partial \left( (\mu_{\bar{f}})_{\bar{f}} \right)} \: . $$
Finally, we note that
$\partial F_{K,\bar{e}}/\partial \mu_{\bar{f}} = \delta_{\bar{e},\bar{f}}$
and $\partial F_K/\partial \mu_{\bar{f}} = 0$ for
$\bar{e}, \bar{f} \subset \partial K$.

%%%%%%%%%%%%%%%%%%%%%%%%%%%%%%%%%%%%%%%%%%%%%%%%%%%%%%%%%%%%%%%%%%%%%%
%
\subsection{Multigrid algorithms for the linearized equations}
\label{subsec:multigrid}
We have shown in Section \ref{sec:equiv} that after elimination 
of the flux variables the hybridized mixed formulation is equivalent
to the nonconforming finite element method (\ref{nonconf.FEM}).
Hence we can apply the multigrid algorithm developed for this type of
nonconforming discretization \cite{Brenner:92} to solve the system 
of linearized equations.

Furthermore, Chen \cite{Chen:96} showed that -- 
at least for linear elliptic problems -- the equations after elimination
of flux and element variables correspond to the equations resulting
from a nonconforming discretization using the Crouzeix--Raviart ansatz space.
Therefore the intergrid transfer operators for the Crouzeix--Raviart ansatz 
space \cite{Braess/Verfuerth:90} can be used to develop a multigrid algorithm
for this linear system.
We do not show this equivalence for the nonlinear problem,
but propose to use this multigrid algorithm.
It showed a good convergence behavior in numerical experiments.

%%%%%%%%%%%%%%%%%%%%%%%%%%%%%%%%%%%%%%%%%%%%%%%%%%%%%%%%%%%%%%%%%%%%%%
%
\section{Application to Darcy--Forchheimer flow in porous media}
\label{sec:dafo}
%
%%%%%%%%%%%%%%%%%%%%%%%%%%%%%%%%%%%%%%%%%%%%%%%%%%%%%%%%%%%%%%%%%%%%%%
The flow of a gas through a porous medium is described 
by the following equations (see e.g.\ \cite{Whitaker:96}):
the Darcy--Forchheimer equation
$$ %\begin{equation}				\label{Da-Fo.eq.orig}
\frac{\mu(\xvec,t)}{k(\xvec)} \, \vvec(\xvec,t)
+ \beta_{\mathrm{Fo}}(\xvec) \, \rho(\xvec,t) \, 
  |\vvec(\xvec,t)| \, \vvec(\xvec,t)
+ \nabla P(\xvec,t) = 0 \: ,
\quad (\xvec,t) \in \Omega \times [0,T] \: ,
$$ %\end{equation}
where $| \cdot |$ denotes the euclidean norm, the continuity equation
$$ %\begin{equation}				\label{cont.eq.orig}
\phi(\xvec) \, \frac{\partial \rho(\xvec,t)}{\partial t} 
+ \Div(\rho(\xvec,t) \vvec(\xvec,t)) = \tilde{f}(\xvec,t) \: ,
\quad (\xvec,t) \in \Omega \times [0,T]
$$ %\end{equation}
and an equation of state, we take the ideal gas law 
$$ %\begin{equation} 				\label{i.g.l}
\rho(\xvec,t) = \frac{P(\xvec,t) \, W(\xvec,t)}{R_0 \, \Theta(\xvec,t)} 
              =: P(\xvec,t) \gamma(\xvec,t)
\: ,
\quad (\xvec,t) \in \Omega \times [0,T] \: .
$$ %\end{equation}
The unknowns here are the pressure $P$, the density $\rho$ 
and the volumetric flow rate $\vvec$ of the gas.
Porosity $\phi$, permeability $k$ 
and Forchheimer coefficient $\beta_{\mathrm{Fo}}$ of the porous medium,
viscosity $\mu$, molecular weight $W$ and temperature $\Theta$ of the gas,
and the universal gas constant $R_0$ are given as well as the source term $\tilde{f}$.
Assuming $\rho>0$ and introducing new variables $p = |P| P$ and 
$\uvec = |\rho| \vvec$ these equations can be transformed into
\begin{subequations}
\begin{align}				\label{Da-Fo.eq}
\left( \alpha(\xvec,t) + \beta(\xvec,t) |\uvec(\xvec,t)| \right) \uvec(\xvec,t)
+ \nabla p(\xvec,t) = & \: 0 \: ,
& (\xvec,t) \in \Omega \times [0,T] \: , \\	\label{cont.eq}
\phi(\xvec) \, \partial_t \rho(p(\xvec,t),\xvec,t) 
+ \Div \uvec(\xvec,t) = & \: \tilde{f}(\xvec,t) \: ,
& (\xvec,t) \in \Omega \times [0,T] \: ,
\end{align}
\end{subequations}
where 
$$ \gamma(\xvec,t) := \frac{W(\xvec,t)}{R_0 \, \Theta(\xvec,t)} \; , \quad
   \alpha(\xvec,t) := \frac{2 \, \mu(\xvec,t)}{\gamma(\xvec,t) \, k(\xvec)} \; , \quad
   \beta(\xvec,t) := \frac{2 \, \beta_{\mathrm{Fo}}(\xvec)}{\gamma(\xvec,t)} \; , $$
and the equation of state $\rho=\rho(p)$ is defined by
\begin{equation}				\label{def:rho(p)}
\rho(p(\xvec,t),\xvec,t) :=
\gamma(\xvec,t) \frac{p(\xvec,t)}{\sqrt{|p(\xvec,t)|}} \: .
\end{equation}
In particular, we consider the mixed finite element discretization
of the semidiscrete problem with Dirichlet boundary conditions.
Discretization in time of the continuity equation (\ref{cont.eq})
with the implicit Euler method yields for the k-th time step
\begin{subequations}				\label{semid.pr}
\begin{align}					\label{semid.Da-Fo.eq}
\left( \alpha^k(\xvec) + \beta^k(\xvec) |\uvec^k(\xvec)| \right) \uvec^k(\xvec)
+ \nabla p^k(\xvec) = & \: 0 \: , & \xvec \in \Omega \: , \\
\frac{\phi(\xvec)}{\Delta t}			\label{semid.cont.eq}
\gamma^k(\xvec) \frac{p^k(\xvec)}{\sqrt{|p^k(\xvec)|}}
- \frac{\phi(\xvec)}{\Delta t} 
  \gamma^{k-1}(\xvec) \frac{p^{k-1}(\xvec)}{\sqrt{|p^{k-1}(\xvec)|}}
+ \Div \uvec^k(\xvec) = & \: \tilde{f}^k(\xvec)
 \: , 
& \xvec \in \Omega \: , \\
p^k(\xvec) = & \: g^k(\xvec) \: , & \xvec \in \partial\Omega \: .
\end{align}
\end{subequations}
The initial condition $p(\cdot,t^0) = p^0(\cdot)$ is given.
Similarly, we use the denotations
$p^k := p(\cdot,t^k)$ and $\uvec^k := \uvec(\cdot,t^k)$ for the unknown 
solutions and, analogously defined, $\alpha^k$, $\beta^k$ and $\gamma^k$ 
for the coefficient functions, $g^k$ for the boundary conditions 
and $\tilde{f}^k$ for the source term.
Note that for each $k \in \{1,\ldots,K\}$ the function $p^{k-1}$ is known.
Hence we can define an augmented source term by
$$ f^k(\xvec) := \tilde{f}^k(\xvec) + \frac{\phi(\xvec)}{\Delta t} 
   \gamma^{k-1}(\xvec) \frac{p^{k-1}(\xvec)}{\sqrt{|p^{k-1}(\xvec)|}} \: . $$

Evidently, (\ref{semid.pr}) is of the form (\ref{mixed.form}).
Omitting the superscript $k$, the nonlinear mappings $R$ and $G$ are defined by
$$ R(p) := \frac{\phi}{\Delta t} \gamma \frac{p}{\sqrt{|p|}}
   \quad \text{and} \quad
   G(\uvec) := \left( \alpha + \beta |\uvec| \right) \uvec \: . $$
We require $\phi, \alpha, \beta, \gamma \in L^\infty(\Omega)$ and additionally 
$$ \left. \begin{array}{l}
0 < \underline{\phi} \le \phi(\xvec) \le \overline{\phi} < \infty \: , \\
0 < \underline{\alpha} \le \alpha(\xvec) \le \overline{\alpha} < \infty \: , \\
0 < \underline{\beta} \le \beta(\xvec) \le \overline{\beta} < \infty \: , \\
0 < \underline{\gamma} \le \gamma(\xvec) \le \overline{\gamma} < \infty
\end{array} \right\}
\mbox{ for almost every } \xvec \in \Omega \: . $$
Then is is easy to show that the mapping
$R: L^{3/2}(\Omega) \to L^3(\Omega)$ is continuous, coercive and 
strictly monotone (cf.\ the proof of Prop.~1.5 b) in \cite{Knabner/Summ}).
Since $\phi$ and $\gamma$ are bounded, 
an application of H\"older's inequality yields
$$ \left\| R(p) - R(q) \right\|_{0,3,\Omega}
  \le \frac{\overline{\phi} \, \overline{\gamma}}{\Delta t} 
   \left\| \frac{p}{\sqrt{|p|}} - \frac{q}{\sqrt{|q|}} \right\|_{0,3,\Omega}
  \le \frac{\overline{\phi} \, \overline{\gamma}}{\Delta t}
   \sqrt{2} \left\| p - q \right\|_{0,3/2,\Omega}^{1/2} $$
for all $p, q \in Q$.
Coercivity and strict monotonicity follow directly from the fact 
that $\phi$ and $\gamma$ are bounded from below by positive constants.
Consequently we have $r=3/2$ and $r'=3$ here.

Similarly, we obtain that
$G: \left( L^3(\Omega) \right)^2 \to \left( L^{3/2}(\Omega) \right)^2$
is a continuous and uniformly monotone mapping 
(see the proof of Prop.~1.2 in \cite{Knabner/Summ}).
Indeed, an application of H\"older's inequality yields
$$ \left\| G \uvec - G \vvec \right\|_{0,3/2,\Omega} \le
   \left( C(\overline{\alpha}) + C(\overline{\beta})
   \left( \|\uvec\|_{0,3,\Omega} + \|\vvec\|_{0,3,\Omega} \right) \right)
   \| \uvec - \vvec \|_{0,3,\Omega} \: , $$
where $C(\overline{\alpha})$ and $C(\overline{\beta})$ are constants depending
only on $\overline{\alpha}$, $\overline{\beta}$ and the domain $\Omega$.
Furthermore,
$$ \int_\Omega (G\uvec - G\vvec) \cdot (\uvec - \vvec) \, d\xvec \ge
   \frac{C(\underline{\beta})}{2} \|\uvec -\vvec\|_{0,3,\Omega}^3 \quad 
   \text{for} \quad \uvec, \vvec \in \left( L^3(\Omega) \right)^2 \: , $$
where $C(\underline{\beta})$ depends only on $\underline{\beta}$ and $\Omega$.
In particular, the uniform monotonicity of $G$ implies that $G$ is strictly
monotone and coercive. Thus $G$ is invertible. 
The inverse mapping $G^{-1}$ is given by
$$ G^{-1}(\vvec) = \frac{\sqrt{\alpha^2 +  4 \beta |\vvec|} - \alpha}
			{2 \beta |\vvec|} \, \vvec \: . $$
Consequently we have $s=3$ and $s'=3/2$ here and the space $V$ is given by
$$ V = W^3(\Div;\Omega) := 
   \left\{ \vvec \in \left( L^3(\Omega) \right)^2 \bigm|
   	   \Div \vvec \in L^3(\Omega) \right\} \: . $$
Some properties of $W^3(\Div;\Omega)$ are established in 
\cite[Appendix]{Knabner/Summ}.

Obviously, $G$ fulfills $G(-\vvec) = -G(\vvec)$ for all $\vvec \in V$.
Furthermore, $G: \R^2 \to \R^2$ is continuously differentiable, 
since all partial derivatives are continuous.
The Jacobian matrix of $G$ is given by 
$DG(\xvec) = (\alpha + \beta |\xvec|) \mathrm{Id} 
	    + \frac{1}{|\xvec|} \xvec \xvec^T$,
where $\mathrm{Id}$ denotes the unit matrix.
In Lemma \ref{lemma:DF:invertible} we employed the assumption
that there is a constant $\underline{a} > 0$ 
and a norm $\| \cdot \|_{V_K}$ on $V_K$ such that 
the operator $A_K: V_K \to V_K'$ defined in \eqref{A_K:def} fulfills
inequality (\ref{ass:DF:invertible}). %, i.e.,
%$$ \left\| A_K \uvec_K - A_K \vvec_K \right\|_{V_K'}
%   \ge \underline{a} \left\| \uvec_K - \vvec_K \right\|_{V_K} 
%   \quad \text{for all } \uvec_K, \vvec_K \in V_K \: . $$
Choosing $\| \cdot \|_{V_K} = \| \cdot \|_{0,2,\Omega}$,
this inequality follows from the observation
\begin{align*}
& \langle A_K \uvec_K - A_K \vvec_K , \uvec_K - \vvec_K \rangle_{V_K' \times V_K} \\
& = \int_K \alpha (\uvec_K-\vvec_K) \cdot (\uvec_K-\vvec_K) \, d\xvec
    + \int_K \beta \left( |\uvec_K| \uvec_K - |\vvec_K| \vvec_K \right) 
             \cdot (\uvec_K-\vvec_K) \, d\xvec \\
& \ge \underline{\alpha} \left\| \uvec_K - \vvec_K \right\|_{V_K}^2 
    \quad \text{for  all } \uvec_K , \vvec_K \in V_K \: ,
\end{align*}
because the second term is non-negative.
In particular, this implies (\ref{ass:DF:invertible}), because
\begin{align*}
\left\| A_K \uvec_K - A_K \vvec_K \right\|_{V_K'}
= & \: \sup_{\wvec_K \in V_K} 
   \frac{\left| \langle A_K \uvec_K - A_K \vvec_K, \wvec_K \rangle \right|}
        {\|\wvec\|_{V_K}} \\
\ge & \: \frac{\left| \langle A_K \uvec_K - A_K \vvec_K , 
		              \uvec_K - \vvec_K \rangle \right|}
              {\left\| \uvec_K - \vvec_K \right\|_{V_K}}
\ge \underline{\alpha} \left\| \uvec_K - \vvec_K \right\|_{V_K} \: .
\end{align*}

For the elimination of element variables, we need the monotonicity 
of the mapping $C_K: \R \to \R$. Here 
$$ C_K(p_K) = \int_K R(p_K) \, d\xvec 
  = \int_K \frac{\phi}{\Delta t} \gamma \frac{p_K}{\sqrt{|p_K|}} \, d\xvec 
  = \int_K \frac{\phi}{\Delta t} \gamma \, d\xvec \, \frac{p_K}{\sqrt{|p_K|}} $$
such that $C_K$ is strictly monotone, since the integral is positive.
Note that $C_K$ is continuously differentiable for all 
$p_K \in \R \setminus \{0\}$ and its derivative 
$C_K'(p_K)=\int_K \phi \gamma \, d\xvec/(2 \, \Delta t \sqrt{|p_K|})$
is not defined for $p_K=0$.
Since inversion of $C_K$ yields
$$ p_K = p_K ( \vec{u}_K )
       = | P_K(\vec{u}_K) | \, P_K (\vec{u}_K) \: , $$
where
$$ P_K (\vec{u}_K) :=
   \frac{\Delta t}{\int_K \phi \gamma \, d\xvec}
   \left( \int_K f \, d\xvec
	 - \sum_{e \subset \partial K} u_{K,e} \right) \, , $$
we can directly compute the partial derivatives 
$$ \frac{\partial p_K}{\partial u_{K,e}} (\vec{u})
   = - 2 \left| P_K (\vec{u}) \right| \,
    \frac{\Delta t}{\int_K \phi \gamma \, d\xvec}
   = - 2 \sqrt{\left| p_K (\vec{u}) \right|} \,
    \frac{\Delta t}{\int_K \phi \gamma \, d\xvec} $$
without using the implicit function theorem.
These partial derivatives are continuous in $\R$.
A closer look into the proof of Lemma \ref{lemma:bar(F)_(K,e):invertible} shows
that we needed the condition on $C_K'$ only to ensure that the partial 
derivatives $\partial p_K/\partial u_{K,e}$ exist and are continuous.
Thus we can replace the condition on $C_K'$ by the requirement that the 
partial derivatives $\partial p_K/\partial u_{K,e}$ exist and are continuous,
which is fulfilled in our application.

\bibliographystyle{amsplain}

\begin{thebibliography}{10}

\bibitem{Arbogast/Chen:95} 
T.~Arbogast and Z.~Chen,
\textit{On the implementation of mixed methods as nonconforming methods for
        second-order elliptic problems},
Math. Comput. \textbf{64} (1995), 943--972.

\bibitem{Braess/Verfuerth:90}
D.~Braess and R.~Verf\"{u}rth,
\textit{Multigrid methods for nonconforming {F}inite {E}lement {M}ethods},
SIAM J. Numer. Anal., \textbf{27} (1990), 979--986.

\bibitem{Brenner:92}
S.~C. Brenner,
\textit{Multigrid algorithm for the lowest-order {R}aviart-{T}homas mixed
        triangular finite element method},
SIAM J. Numer. Anal., \textbf{29} (1992), 647--678.

\bibitem{Brezzi/Fortin}
F.~Brezzi and M.~Fortin,
\textit{Mixed and hybrid finite element methods}, 
Springer series in computational mathematics, vol.~15,
Springer-Verlag, Berlin, Heidelberg, New York, 1991.

\bibitem{Chen:96}
Z.~Chen,
\textit{Equivalence between and multigrid algorithms for nonconforming and
        mixed methods for second-order elliptic problems},
East-West J. Numer. Math., \textbf{4} (1996), 1--33.

\bibitem{Knabner/Summ}
P.~Knabner and G.~Summ,
\textit{Solvability of the mixed formulation for Darcy--Forchheimer flow in porous media},
submitted to M2AN Math. Model. Numer. Anal.

\bibitem{Whitaker:96}
S.~Whitaker,
\textit{The {F}orchheimer equation: A theoretical development},
Transp. Porous Media, \textbf{25} (1996), 27--61.

\bibitem{Zeidler}
E.~Zeidler
\textit{Nonlinear functional analysis and its applications --
        Nonlinear monotone operators},
Springer-Verlag, Berlin, Heidelberg, New York, 1990.
\end{thebibliography}

\end{document}